\documentclass[letterpaper, 10 pt, onecolumn]{ieeeconf}

\IEEEoverridecommandlockouts
\usepackage{graphicx}
\usepackage{amsmath}
\usepackage{amssymb}
\usepackage{amsfonts}
\usepackage{epsfig}
\usepackage{color}

\title{\LARGE \bf Consensus in non-commutative spaces}

\author{Rodolphe Sepulchre, Alain Sarlette and Pierre Rouchon\thanks{R.~Sepulchre and A.~Sarlette are with the Department of Electrical Engineering and Computer Science (Montefiore Institute, B28), University of Li\`ege, 4000 Li\`ege, Belgium. P.~Rouchon is with the Centre Automatique et Syst\`emes, Mines ParisTech, 75272 Paris cedex 06, France. {\tt\footnotesize r.sepulchre@ulg.ac.be, alain.sarlette@ulg.ac.be, pierre.rouchon@mines-paristech.fr}}%
}

\begin{document}

\maketitle
\thispagestyle{empty}
\pagestyle{empty}


\begin{abstract}
	
Convergence analysis of consensus algorithms is revisited in the light of the Hilbert distance.  Tsitsiklis Lyapunov function is shown to be the Hilbert distance to consensus in log coordinates. Birkhoff theorem, which proves contraction of the Hilbert metric for any positive homogeneous monotone map, provides an early yet general convergence result for consensus algorithms. Because Birkhoff theorem holds in arbitrary cones, we extend consensus algorithms to the cone of positive definite matrices. The proposed generalization finds applications in the convergence analysis of quantum stochastic maps, which are a generalization of stochastic maps to non-commutative probability spaces.

\end{abstract}


\section{Introduction}

The convergence of consensus algorithms in linear spaces has attracted much interest in recent years due to its many applications in distributed computation and control. Because the matrix defining the update of the algorithm can be given a probabilistic interpretation as a stochastic map or a graph interpretation as an adjacency matrix, the topic has fostered a number of interesting connections between control theory, probability theory, and graph theory. Our paper revisits some of those connections in the light of a classical fixed point result and uses this link to define a generalization of  consensus algorithms in non-commutative spaces.

A key component of our approach is a 1957 theorem of Garrett Birkhoff. Birkhoff theorem provides a general fixed point result for homogeneous monotone positive maps defined on closed cones of Banach spaces. The key idea is to use a metric
introduced by Hilbert as contraction measure. Hilbert's metric is projective, that is, it measures distances between homogeneous rays of the cone.

In the positive orthant, Hilbert distance between vectors $x$ and $y$ is  $d_H(x,y)= \max \log (x_i/y_i) - \min \log (x_i/y_i)$,
which is indeed invariant by scaling of either $x$ or $y$. Taking $y = \mathbf{1}= [1, \dots, 1 ]^T$, this provides a natural distance to consensus; in probabilistic terms this corresponds to a uniform density vector. In fact, Hilbert distance to consensus coincides (in $log$ coordinates) with the  Lyapunov function introduced by Tsitsiklis to study the convergence of consensus algorithms. Row-stochastic matrices define linear positive maps in the positive orthant. As a consequence,  Birkhoff theorem provides an early convergence result for  consensus algorithms.

Because Birkhoff theorem applies in arbitrary cones, it opens the way to several generalizations of consensus theory. The present paper focuses on the cone of hermitian positive definite matrices. Our motivation stems from the generalization of probability theory in the quantum (non-commutative) setting. While stochastic maps operate on probability \emph{vectors} (elements of the positive orthant), quantum stochastic maps operate on probability \emph{matrices} (elements of the cone of hermitian positive semidefinite matrices).  Quantum stochastic maps are unital maps in this cone, i.e. they map identity on the identity. Their dual  (Kraus maps) are trace-preserving and  appear frequently in the study of  open quantum systems.

In the cone of hermitian positive definite matrices, Hilbert distance from a matrix $X$ to the identity matrix is $d_H(X,I)=  \log (\lambda_{\max}) -
\log (\lambda_{\min})$. As in the positive orthant, this distance is projective, i.e. invariant by scaling (of either $X$ or $I$). For unital maps, it provides a natural distance to
the uniform probability density matrix $I$, which defines ``consensus'' in this context. Because quantum stochastic maps define positive linear monotone maps on the cone of positive definite matrices, Birkhoff theorem shows that the distance to consensus can only decrease when measured with the Hilbert metric. The paper explores some implications of this basic property in  the context of quantum channels. The preliminary examples discussed in this paper indicate that the analogy between the classical consensus theory in commutative spaces and the proposed consensus theory in non-commutative spaces offers new avenues to characterize the convergence properties of quantum channels. In particular, Tsitsiklis approach to consensus theory proved especially useful in asymmetric and/or time-varying settings where a more conventional approach based on quadratic Lyapunov functions fails. Likewise, the Hilbert distance as a measure of distance to consensus in a quantum setting could prove particularly useful in the context of  asymmetric and/or time-varying quantum maps.

The paper is organized as follows: the convergence theory of consensus algorithms is briefly revisited in Section 2. Birkhoff theorem is recalled in Section 3. Its interpretation in the framework of consensus algorithms is discussed in Section 4. Section 5 extends the framework to the cone of positive definite matrices. The notion of non-commutative consensus, motivated by quantum applications, is introduced in Section 6. Section 7 contains concluding remarks and future perspectives.


\section{Consensus algorithms and Tsitsiklis Lyapunov function}

Linear consensus algorithms give rise to time-varying systems
\begin{equation}\label{Syst1}
x(t+1) = A(t) x(t),  \;  \; x(t) \in \mathbb{R}^n
\end{equation}
where for each $t= 0, 1, \dots$ the matrix $A(t)$ is row-stochastic, i.e. its elements are non-negative ($a_{ij}(t) \ge 0$), and it preserves the unit element $\mathbf{1}=(1, \dots, 1)^T$ ($A(t) \mathbf{1}=  \mathbf{1}$). The graph interpretation is that $N$ nodes exchange information about a scalar quantity $x_i(t)$ along communication edges weighted by the coefficients $a_{ij}(t)$. ($a_{ij}(t)=0$ means that no information is transferred from node $j$ to node $i$ at time $t$). The update of a given node is a weighted average of its own value and the values communicated by neighbors:
\begin{equation}\label{update}
x_i(t+1) = \sum_{j=1}^n a_{ij}(t) x_j(t),  \; \; \sum_{j=1}^n a_{ij}(t)=1
\end{equation}

Remarkably, uniform convergence of the time-varying system (\ref{Syst1}) can be characterized quite precisely. The early analysis of Tsitsiklis \cite{tsitsiklis86} rests on the basic but  fundamental observation that the (time-invariant) Lyapunov function
\begin{equation}\label{tsitsiklisLyap}
V(x) = \max_{1 \le i \le n}{x_i} - \min_{1 \le i \le n}{x_i}
\end{equation}
is non-increasing along the solutions of (\ref{Syst1}). The proof that it decreases {\it strictly} over a uniform horizon under appropriate assumptions  only involves elementary calculations (see e.g. \cite{blondel05} for details).

An example of general convergence result for (\ref{Syst1})  is as follows:

\noindent \textbf{Theorem 1:} \emph{Let $A(t), \; t \ge 0$  be a sequence of $ n \times n $ row-stochastic matrices and
suppose that the following conditions hold:
\begin{itemize}
\item all nonzero $a_{ij}(t)$ are larger than a certain $\alpha > 0$.
\item All diagonal elements satisfy $a_{ii}(t)>0$
\end{itemize}
Then the solution of (\ref{Syst1}) converges
exponentially to an equilibrium co-linear to  $\mathbf{1} $ if and only if
there exists a finite horizon $ T > 0 $ such that for all $t_0 \in  \mathbb{N}$,  there exists a node connected to
all other nodes across $[t_0, t_0 + T]$.
}
\vspace{2mm}

The proof of (variants) of this result has appeared in various papers, including important contributions by Tsitsiklis \cite{tsitsiklis86}, Jadbabaie et al. \cite{jadbabaie03}, and Moreau \cite{moreau05}.

It is also remarkable that the convergence analysis of the linear system (\ref{Syst1}) cannot be established by means of a time-invariant {\it quadratic} Lyapunov function. The recent paper \cite{olshevsky08} provides the explicit construction of a sequence of eight matrices $A(t)$ that satisfies the assumptions of Theorem 1 but that does not admit a common time-invariant quadratic Lyapunov function. In other words, the non-quadratic nature of the Lyapunov function (\ref{tsitsiklisLyap}) is essential and a distinctive feature of the analysis of consensus algorithms.


\section{Birkhoff Theorem}

As will be shown in the next section, the Lyapunov function (\ref{tsitsiklisLyap}) is closely related to  Hilbert metric. This connection follows from a result of Birkhoff \cite{birkhoff57} that  we now summarize, following the exposition in \cite{bushell73}.

Let ${\cal X}$ be a real Banach space and let ${\cal K}$ be a closed solid cone in ${\cal X}$, that is, a closed subset ${\cal K}$ with the properties: (i) ${\rm int}{\cal K}$, the interior of ${\cal K}$, is not empty, (ii) ${\cal K}+{\cal K} \subset {\cal K}$, (iii) $\lambda {\cal K} \subset {\cal K}$ for all $\lambda \ge 0$, and (iv) ${\cal K} \cap -{\cal K} = \{0\}$. The partial order $x \preceq y$ means that $y-x \in {\cal K}$.

For $x, y \in {\cal K}_0 = {\cal K} \backslash \{0 \}$, define $M(x,y) = \inf \{ \lambda:  x - \lambda  y \preceq 0 \} $
and $
 m(x,y) = \sup \{ \lambda: x - \lambda y  \succeq 0 \} $. Then Hilbert metric $d(\cdot,\cdot)$ is defined in ${\cal K}_0$ by
 \begin{equation}\label{Hilbert}
 d(x,y) = \log \{M(x,y)/m(x,y) \}
\end{equation}
Because of the invariance property $d(\alpha x, \beta y)=d(x,y)$ for all $\alpha, \beta >0$, Hilbert metric is, strictly speaking, a
{\it projective} metric, i.e. a distance between the equivalence classes (rays) $[x]$ and $[y]$, where
 $[x] =  \{ \alpha x \mid  \alpha >0 \}$.

 Birkhoff theorem is a contraction result for positive maps in ${\cal K}$. A map $A: {\cal K} \rightarrow {\cal K}$ is said to be {\it non-negative}; a map
 $A: {\rm int}{\cal K} \rightarrow {\rm int}{\cal K}$ is said to be {\it positive}. If $A$ is positive, the projective diameter $\Delta(A)$ of $A$ is defined by
 $$ \Delta(A) = \sup \{ d(Ax,Ay): x, y \in {\rm int}{\cal K} \} \; .$$

 \noindent \textbf{Theorem 2 [Birkhoff, 1957]:} \emph{Let $A$ a map in ${\cal K}$ that satisfies the following properties:}
 \begin{itemize}
 \item[(i)] \emph{$A$ is positive, i.e. $A: {\rm int}{\cal K} \rightarrow {\rm int}{\cal K}$.}
 \item[(ii)] \emph{$A$ is homogeneous of degree $p$ in ${\rm int}{\cal K}$, i.e. $A(\lambda x) = \lambda^p A(x)$ for all $\lambda >0$.}
 \item[(iii)] \emph{$A$ is monotone, i.e. $ x \preceq y  \Rightarrow A(x) \preceq A(y) $.}
 \end{itemize}
 \emph{Then $A$ is a contraction for the Hilbert metric and the following holds:
 $ \forall x,y \in {\rm int}{\cal K} :  \newline d(A(x),A(y)) \le \tanh (\frac{1}{4} \Delta(A)) \; d(x,y) $ }
\vspace{2mm}

It should be noted that Banach contraction map theorem applied to Birkhoff theorem provides a far-reaching
generalization of Perron-Frobenius theorem: if $\Delta (A) < \infty$ and if the metric projective space $E= ({\rm int}{\cal K}  \backslash \tilde \; , d)$ is
complete, then there exists a unique eigenvector of $A$ in $E$.

Hilbert metric is only one among the many metrics contracted by positive monotone maps. Another such (closely related) metric is
Thompson metric $d_T(x,y)=  \log \max \{ M(x,y), m^{-1}(x,y) \} $. However, Hilbert metric typically provides the best contraction ratio, see \cite{bushell73} for details.


\section{Birkhoff Theorem in the positive orthant and Consensus algorithms}

The most direct application of Birkhoff theorem is in the positive orthant. In this case, ${\cal X}= \mathbb{R}^n$ and $ {\cal K} = \{(x_1, \dots, x_n): x_i \ge 0,\; 1 \le i \le n \} $.
Then $M(x,y)= \max_i (x_i/y_i)$ and $m(x,y) = \min_i (x_i/y_i)$ and hence Hilbert metric expresses as
$$ d(x,y) = \log \frac{\max_i(x_i/y_i)}{\min_i(x_i/y_i)} \; .$$

By definition, stochastic matrices $A(t)$ define positive monotone maps in ${\cal K}$. By linearity, they are also homogeneous of degree one.

Because $\mathbf{1}$ is a fixed point for all $A(t)$, Birkhoff theorem implies that
$d(x,\mathbf{1}) \le  \tanh (\frac{1}{4} \Delta( A(t))) \; d( A(t) x,\mathbf{1})$ for all $t$. This means
that for each $x \in {\rm int}{\cal K}$, Birkhoff theorem provides the Lyapunov function
\begin{equation}\label{birkhoffLyap}
 V_B(x) = d(x,\mathbf{1}) = \log \frac{\max_i(x_i)}{\min_i(x_i)} = \max \log (x_i) - \min \log (x_i)
\end{equation}
This is exactly the Lyapunov function (\ref{tsitsiklisLyap}) in $\log$ coordinates. Both Lyapunov functions provide a measure of the
projective distance from $[x]$ to $[\mathbf{1}]$. Tsitsiklis Lyapunov function is   translation-invariant  ($V(x+ \lambda \mathbf{1}) = V(x)$) while
Birkhoff Lyapunov function is scaling-invariant ($V_B(\lambda x)=V(x)$ for $\lambda>0$). They both provide a measure of the diameter
of the convex hull of $(x_1, \dots, x_n)$, which is the Lyapunov function proposed by Moreau \cite{moreau05}.

Birkhoff theorem is thus an early contribution that points to two important facts: (i) the natural distance to study the convergence of consensus algorihms is not quadratic and (ii) the convergence analysis should concentrate on the contraction of the interval $[x_{\min},x_{\max}]$.

Birkhoff theorem also provides a contraction coefficient. It is given by  the diameter of $A(t)$
$$  \Delta(A(t)) = \sup \left\{ \log(\frac{ a_{ij}(t) a_{pq}(t)}{ a_{iq}(t) a_{pj}(t)}) \; : 1 \le i,j,p,q \le n \right\}\; . $$
When the diameter is finite, the contraction coefficient is strictly smaller than one. Uniform convergence only requires contraction over a uniform time-horizon. As a consequence,  if $A(t)$ does not have a finite diameter at each instant, it is sufficient for uniform convergence that the finite product $\tilde A(t)=A(t+T)\dots A(t+1) A(t)$ has a finite diameter.

In the symmetric case $A(t)=A^T(t)$, the existence of a finite horizon over which the diameter is finite is closely linked to the existence of a finite horizon over which the communication graph is connected. This link is studied in detail in \cite{jadbabaie03}. It provides a close connection between Birkhoff theorem and Theorem 1 in the {\it symmetric} setting. Symmetric stochastic matrices are called doubly-stochastic matrices (they are row-stochastic and column-stochastic) and symmetric graphs appear most notably when the edges are undirected.

In contrast, the {\it asymmetric} setting (that occurs in directed communication graphs) includes situations where, despite a uniform convergence to consensus, the diameter is infinite over any finite horizon. In those situations, illustrated in the next example, Birkhoff theorem does not provide a strict contraction measure even though elementary explicit calculations show the contraction of $[x_{\min},x_{\max}]$ under the
assumptions of Theorem 1.

\noindent \textbf{Example 1: A consensus map with infinite diameter.}  The matrix
$$A = \left(\begin{array}{cc} 1 & 0 \\ \gamma^2  & 1-\gamma^2 \end{array} \right) \qquad \text{with } \gamma \in (0,1) \; .$$
satisfies the assumptions of Theorem 1.  For an initial condition $(x_1,x_2)$ in the positive orthant, the asymptotic (consensus) value is $(x_1,x_1)$. Nevertheless,  the diameter of $A$ is infinite. Because $A$ is lower-triangular, any finite power $A^k$ also has an infinite diameter. $\hfill \square$

Quantifying the convergence rate with other coefficients than Birkhoff diameter remains a topic of research, see e.g. \cite{AngeliBliman08} and references therein.

It should be emphasized that Birkhoff theorem  is not restriced to linear maps (see \cite{cortes08,moreau05} for examples of nonlinear consensus algorithms).


\section{Birkhoff Theorem in the cone of positive definite matrices}

We now turn to the application of  Birkhoff theorem  in the  cone of positive definite matrices.
In this case, ${\cal X}= \{X=X^* \in \mathbb{C}^{n\times n} \}$ is the set of hermitian matrices, where $\phantom{A}^*$ denotes complex conjugate transpose,
 and $ {\cal K} = \{X \succeq 0 \mid X \in {\cal X} \} $.

Then for $X ,Y \succ 0$,  $M(X,Y)= \lambda_{\max} (XY^{-1})  $ and $m(x,y) = \lambda_{\min} (XY^{-1})$.  Hence Hilbert metric expresses as
$$ d(X,Y) = \log \frac{\lambda_{\max} (XY^{-1})}{\lambda_{\min} (XY^{-1})}.$$
The reader will note that an equivalent expression is
\begin{equation}
\label{hilbertSDP}
 d(X,Y) = \log \frac{\lambda_{\max} (Y^{-1/2}XY^{-1/2})}{\lambda_{\min} (Y^{-1/2}XY^{-1/2})}
 \end{equation}
which shows the close connection of Hilbert metric to the natural invariant metric
\begin{eqnarray}
\label{rieman}
d_R(X,Y) & = & \parallel \log  (Y^{-1/2}XY^{-1/2}) \parallel_F
= {\textstyle \sum_{i=1}^n} \;  \log \lambda_i(Y^{-1/2}XY^{-1/2}) \; .
\end{eqnarray}
In contrast to the Hilbert metric, the metric (\ref{rieman}) is a Riemannian metric. Both metrics share congruence invariance properties: $d (FXF^*,FYF^*)=d (X,Y)
$ for any invertible linear map $F$ in $\mathbb{C}^{n\times n} $.
The Riemannian metric coincides with the Fisher-information metric in statistics \cite{smith05}
 and with the  self-concordant logarithmically
homogeneous barrier ($-{\rm log}\, {\rm det}\, A$)  in optimization \cite{nesterov02}. It makes the cone of positive definite matrices geodesically complete.

Birkhoff theorem implies that any homogeneous monotone positive map $A$ on the cone of positive definite matrices contracts the metric (\ref{hilbertSDP}).


\section{Non-commutative consensus}


\subsection{Generalizing stochastic maps to non-commutative spaces}

A key feature of consensus algorithm (\ref{Syst1}) on $\mathbb{R}^n$ is the stochasticity of maps $A(t)$, i.e. (a) they have non-negative elements and (b) $A(t) \mathbf{1}=  \mathbf{1}$. An interesting  analog of stochastic maps on the cone ${\cal K}$ of positive semi-definite  hermitian matrices of order $n$ is  given by the dual maps of  Kraus maps describing quantum channels (see \cite{NielsenChuang} chapter 8 for a tutorial presentation of quantum channels; see \cite{HarocheBook} chapter~4 for a physical presentation of decoherence and  Kraus maps).
These dual maps are \emph{unital completely positive maps}, which admit the following characterization \cite{Choi75}: a linear map $\Phi $ on ${\cal K}$ is unital completely positive if and only if it admits the expression
  \begin{eqnarray}
\label{kraus1}
\forall X\in{\cal K},\;\;\;\quad ~\Phi(X) & = & \sum_i V_i^* X V_i
 \qquad \text{with} \; \; \sum_i V_i^* V_i = I
 \end{eqnarray}
for matrices $V_i \in \mathbb{C}^{n \times n}$. Completely positive maps, characterized by the first line of (\ref{kraus1}), map ${\cal K}$ into itself. They are positive linear maps on ${\cal K}$ in the terminology of the present paper.   Indeed, $X \succ 0$ implies that, for $w\in\mathbb{C}^{n}$,  $w^* \Phi(X) w = \sum_i \; (V_i w)^* X (V_i w) = 0$ if and only if $V_i w = 0$ $\forall i$. But then $\sum_i \; (V_i w)^* V_i w = 0$ and thus $w=0$ since $\sum_i \; V_i^* V_i = I$. Complete positivity of $\Phi$ generalizes the non-negative elements property of stochastic matrices. The additional  unital property $\Phi(I) = I$ is the analog of $A(t)\mathbf{1}= \mathbf{1}$ for stochastic matrices.

This suggests the following non-commutative extension of consensus algorithms:
 \begin{equation}\label{Syst2}
X(t+1) = \Phi_t (X(t)),  \;  \; X(0)=X_0 \in {{\rm int}\cal K}
\end{equation}
where for each $t$, the map $\Phi_t$ is of the form (\ref{kraus1}). The linear consensus algorithm on $\mathbb{R}^n$ can be obtained as a special case of (\ref{Syst2}), see Section \ref{sec:sec:equalsR}.

Note that, in the same way as algorithm (\ref{Syst1}) can be extended from the positive orthant to the entire space   $\mathbb{R}^n$, the non-commutative generalization (\ref{Syst2}) can be extended to arbitary hermitian matrices in ${\cal X} = \{X=X^* \in \mathbb{C}^{n\times n} \}$. This is a consequence of the translation invariance property of the algorithm: the property $\Phi_t(I) = I$ and linearity imply that, if $X(t)$ is the solution of (\ref{Syst2}) with $X(0)=X_0 \in {\cal K}$, then the solution of (\ref{Syst2}) with initial condition $X_0 + \alpha \, I$ equals $X(t) + \alpha \, I$, for any $\alpha \in \mathbb{R}$ (in particular $\alpha < 0$).


\subsection{Birkhoff theorem and convergence of non-commutative consensus}

By definition, $\Phi_t$ is a linear monotone positive map on the cone of positive definite matrices (i.e. satisfying conditions (i)-(iii) of Theorem 2). The following generalizes  to the non-commutative setting (\ref{Syst2}) the basic convexity property, $\{x_i(t+1), \; i=1,2...n \} \subset [\max(x_i(t)),\, \min(x_i(t))]$ for solutions of (\ref{Syst1}).

\noindent \textbf{Lemma 1:} \emph{Let $\Phi$ be a map of the form (\ref{kraus1}). Then $\lambda_{\min}(\Phi(X)) \geq \lambda_{\min}(X)$ and $\lambda_{\max}(\Phi(X)) \leq \lambda_{\max}(X)$} for any $X=X^*\in \mathbb{C}^{n\times n}$.

\underline{Proof:}
Define $Y_1 := X-(\lambda_{\min}(X)-\varepsilon)\, I$, for any $0 < \varepsilon \ll 1$. Then $Y_1 \in {\rm int}{\cal K}$ implies $\Phi(Y_1) \in {\rm int}{\cal K}$. Since $\Phi(Y_1) := \Phi(X) - (\lambda_{\min}(X)-\varepsilon)\, I$ we have the first bound when $\varepsilon \rightarrow 0$. The second bound is obtained analogously with $Y_2 := -X + (\lambda_{\max}(X)+\varepsilon)\, I$. $\hfill \square$

The following result is a direct consequence of Birkhoff theorem and applying Lemma 1 recursively at each $t$.\vspace{2mm}

\noindent \textbf{Theorem 3:} \emph{Let $\Phi_t, \; t \ge 0$  be a sequence of maps of the form (\ref{kraus1}). Then the Lyapunov function} $$ V(X) = d(X, I)= \log \frac{\lambda_{\max}(X)}{\lambda_{\min}(X)} $$
\emph{is non-increasing along the solutions of (\ref{Syst2}). Denoting $\tilde \Phi_{[t,t+T]} := \Phi_{t+T} \circ \dots \circ \Phi_{t+1} \circ \Phi_t$, the Lyapunov function  satisfies
$$ V(X(t+T)) \le  \tanh(\tfrac{1}{4} \Delta(\tilde \Phi_{[t,t+T]}) ) \; V(X) $$
In particular, if the maps $\tilde \Phi_{[t,t+T]}$ have a finite projective diameter  for some finite horizon $T$, then the decrease of $V$ is strict over $T$ iterations of (\ref{Syst2}),  and the solution of (\ref{Syst2}) uniformly converges to a point in the set $\{\lambda I \mid \lambda \in [\lambda_{\min}(X_0),\lambda_{\max}(X_0)] \}$.}\vspace{3mm}

The following proposition is useful to characterize the diameter of the application. \vspace{2mm}

\noindent \textbf{Proposition 1:} \emph{Let $\tilde \Phi$ an abbreviated notation for the map $\tilde \Phi_{[t,t+T]}$ corresponding to a finite number of iterates of (\ref{Syst2}) as in Theorem 3.}
\newline (a) \emph{Define $R(\tilde \Phi) = {\rm sup}\{\log \tfrac{\lambda_{\max}(\tilde \Phi(X))}{\lambda_{\min}(\tilde \Phi(X))} \mid X \in {\rm int}{\cal K} \}$. Then
$R(\tilde \Phi) \leq \Delta(\tilde \Phi) \leq 2 \, R(\tilde \Phi)$.}
\newline (b) \emph{If $\tilde \Phi(P)$ has full rank for all rank-1 projectors $P \in {\cal P}=\{ x x^* \mid x \in \mathbb{C}^n,\; \Vert x \Vert_2 = 1 \}$, then the superior value in $R(\tilde \Phi)$ is obtained for $X \in {\cal P}$ and $R(\tilde \Phi) = \bar{R} < \infty$.}

\underline{Proof:} (a) $\Delta(\tilde \Phi) \geq R(\tilde \Phi)$ follows by choosing $Y=I$ in the definition of $\Delta$. $\Delta(\tilde \Phi) \leq 2 \, R(\tilde \Phi)$ follows from the triangle inequality on the Hilbert metric.

(b) If $\tilde \Phi(P)$ has full rank for all $P \in {\cal P}$, then $R_P(\tilde \Phi) := {\rm sup}\{\log \tfrac{\lambda_{\max}(\tilde \Phi(P))}{\lambda_{\min}(\tilde \Phi(P))} \mid P \in {\rm int}{\cal P} \}$ takes a finite value, since ${\cal P}$ is compact and $P\mapsto \log \tfrac{\lambda_{\max}(\tilde \Phi(P))}{\lambda_{\min}(\tilde \Phi(P))}$ is continuous. Every $X \in {\cal K}$ can be expressed as a sum of positively weighted rank-1 projectors. Since $\lambda_{\max}(X)$ and $\lambda_{\min}(X)$ are respectively convex and concave functions of $X \in {\cal K}$ and $d(X,I)$ is invariant by scaling of $X$, we directly have $R_P(\tilde \Phi) = R(\tilde \Phi)$. $\hfill \square$
\vspace{2mm}

Note that a finite diameter for the finitely iterated map is only a sufficient condition to prove strict contraction. Very much as in the classical consensus setting, there are situations where the strict decay of the Lyapunov function will be established
by direct computations despite an infinite diameter over any finite horizon.


\subsection{Consensus in $\mathbb{R}^n$ as a special case of non-commutative consensus}\label{sec:sec:equalsR}

There is an obvious bijection between the positive orthant of $\mathbb{R}^n$ and the set of \emph{diagonal} positive semidefinite matrices ${\cal D} = \{ X \in \mathbb{R}^{n\times n} \mid X \succeq 0, \, x_{ij} = 0 \text{ if } i\neq j  \}$. In this section we show how consensus in $\mathbb{R}^n$ can be recovered as a special case of the non-commutative generalization, when the state is restricted to diagonal matrices in some fixed basis. These diagonal matrices do all \emph{commute} with each other; this justifies the appellation ``non-commutative'' consensus for the generalized case of non-diagonal matrices.

Denote by ${\cal S}$ the set of all permutation matrices, i.e. matrices with exactly one element equal to $1$ on each row and column, and all other entries $0$. Note that $S^* S = I$ for any $S \in {\cal S}$. Hence a map
\begin{eqnarray}
\label{eq:eqR} X & \rightarrow & \Phi_{\mathcal{S}}(X) = \sum_i \, W_i \, S_i^* X S_i W_i
 \quad \text{with } S_i \in {\cal S} \, , \; W_i \in {\cal D} \; \forall i
\quad \text{and } \sum_i  W_i^2 = I
\end{eqnarray}
is of the form (\ref{kraus1}). Moreover, $\Phi_{\mathcal{S}}$ defined by (\ref{eq:eqR}) maps ${\cal D}$ into itself. Indeed, the map $X \rightarrow S_i^* X S_i$ permutes in the same way both the rows and columns of $X$; for $X \in {\cal D}$ this amounts to permuting its diagonal elements. The $W_i$ further scale each element. This suggests a way to write the classical linear consensus (\ref{Syst1}) in the form of (\ref{Syst2}) on ${\cal D}$ where $x_i$ of (\ref{Syst1}) is identified with $x_{ii}$ of (\ref{Syst2}).  For instance, take $S_i$ to be the permutation sending $x_{jj}$ to position $(j\oplus i,j\oplus i)$ for all $j$, where $\oplus$ denotes addition modulo $n$. Further take element $(j\oplus i,j\oplus i)$ of $W_i$ equal to $\sqrt{a_{j\oplus i,\, i}}$. Then (\ref{eq:eqR}) with this setting is equal to (\ref{Syst1}); in particular, $\sum_i \; (W_i^2)_{jj} = \sum_i \; a_{ji} = 1$. (This  construction is not unique.)


\subsection{Birkhoff theorem on dual and convergence of iterated Kraus maps}

For simplicity, we consider a fixed (time-invariant) map $\Phi$ for the rest of the paper. The adaptation to the time-varying case is straightforward but  left for an extended version of the paper.

The dual map $\Psi$  of $\Phi$ defined by~\eqref{kraus1} is the usual framework to describe the evolution of the density matrix characterizing a state in a quantum channel, see e.g. \cite{NielsenChuang}. It is called a Kraus map and it admits the  form
\begin{eqnarray}
\label{kraus}
\forall Z\in{\cal K},\;\;\;\quad ~\Psi(Z) & = & \sum_i V_i Z V_i^*
\end{eqnarray}
with still $ \sum_i V_i^* V_i = I$.   The \emph{unital} property of $\Phi$ becomes a \emph{trace-preserving} property for the dual map $\Psi$. Kraus maps are thus  \emph{trace preserving completely positive maps}. In quantum applications, the density matrix typically evolves according to
 \begin{equation}\label{QChannel}
Z(t+1) = \Psi(Z(t)),  \;  \; Z(t) \in {\cal K}_1  := \{ X \in {\cal K} \mid  {\rm tr}(X)=1 \}
\end{equation}
where  the $V_i$ defining $\Psi$ are  matrices in~${\mathbb C}^{n\times n}$. From  Theorem~3 and Proposition~1, we can deduce the following convergence result for $\Psi$.
\vspace{3mm}

\noindent \textbf{Theorem 4:} \emph{Let $\Psi$  defined by~(\ref{kraus}). Assume that $N$-iteration of  its dual map $\Phi$ defined by~(\ref{kraus1}) admits a finite projective diameter} $\Delta(\Phi^N) < \infty$, \emph{for some integer $N$}. \emph{ Then, there exists a unique fixed point $\bar Z$ of $\Psi$ in ${\cal K}_1$, and for any $Z_0 \in {\cal K}_1$, we have $\lim_{t\rightarrow +\infty} Z(t)=\bar Z$. For any initial $X_0 \in {\cal X}$ of $X(t)$ defined by~\eqref{Syst2}, we have
$\lim_{t\rightarrow +\infty} X(t) =  {\rm tr}(\bar Z  X_0) I$. Moreover the convergence is exponential.}\vspace{3mm}

\underline{Proof:}
$\tilde \Phi_t$ satisfies assumptions of Theorem~3. Thus for any initial hermitian matrix $X_0$, $X(t)$ converges exponentially towards $f(X_0) I$ where $f(X_0) \in [\lambda_{\min}(X_0),\lambda_{\max}(X_0)]$ is the consensus value, which depends on $X_0$. Since $\tilde \Phi_t$ is a linear map on ${\cal X}$, $f$ is also a linear map on ${\cal X}$. Therefore there exists a unique hermitian matrix  $\bar Z$ such that $f(X_0)={\rm tr}(\bar Z X_0)$ $\forall X_0 \in {\cal X}$. For all $t$ we have, by duality,  ${\rm tr}(Z(t)\, X_0) = {\rm tr}(Z_0\, X(t))$. Thus for any $X_0 \in {\cal X}$,
 ${\rm tr}(Z(t)\, X_0)$ converges exponentially towards ${\rm tr}( Z_0 \, f(X_0) I) = f(X_0) \, {\rm tr}(Z_0) =f(X_0)={\rm tr}(\bar Z X_0)$. Thus $Z(t)$ converges exponentially towards $\bar Z$, which is positive semidefinite since ${\cal K}$  is closed. $\text{tr} (\bar Z)=1$  since $\text{tr}(Z(t))\equiv 1$. $\hfill \square$\vspace{2mm}

The relation between primal and dual maps can also be made for the classical consensus in $\mathbb{R}^n$. Indeed, viewing (\ref{Syst1}) as the dual, the primal is
\begin{equation}\label{eq:ConsDual} z(t+1) = A^* z(t),  \;  \; z(t) \in \mathbb{R}^n \end{equation}
where $A^*$ is now column stochastic. This implies that $\sum_i \; z_i(t)$ is conserved. Analogously to Theorem 4, the scalar $z^* A x$ can be expressed as evolving either on the primal or on the dual; thus if (\ref{Syst1}) converges to consensus for all $x_0$, then (\ref{eq:ConsDual}) converges to some $\bar{z}$ (not necessarily consensus).


\subsection{Illustration on quantum channel applications}

\noindent \textbf{Example 2: spin rotations.} We consider an open quantum channel where a stochastic environment influences a system with density matrix $\rho \in \mathbb{C}^{2 \times 2}$ --- that is a two-level system like e.g. a spin or qbit (see \cite{HarocheBook,NielsenChuang}). Generalizing \cite{NielsenChuang}, we assume that at each iteration, two different ``spin-rotations'' can be applied to the system with probabilities $p$ and $1-p$ depending on the stochastic environment. The expectation $Z \in \mathbb{C}^{2 \times 2}$ of all stochastic (Monte-Carlo) trajectories for $\rho$ is then characterized by
\begin{eqnarray}
\label{eq:spinrot}	&Z(t+1)  =  V_0 Z V_0^* + V_1 Z V_1^* &\\  &\text{with}\quad
V_0 :=  r_z = \sqrt{p} \left( \begin{array}{cc} e^{i \alpha} & 0 \\ 0 & e^{-i \alpha}
	\end{array}\right),
\quad
V_1 :=  r_x = \sqrt{1-p} \left( \begin{array}{cc} \cos(\beta) & i \, \sin(\beta) \\ i \, \sin(\beta) & \cos(\beta)
	\end{array}\right) \, .&
\end{eqnarray}
Note that $\sum_k \; V_k^* V_k = \sum_k \; V_k \, V_k^* = I$. This implies that the quantum stochastic map is both trace-preserving \emph{and} unital. In fact the Kraus map $\Phi$ and its dual $\Psi$ are equal, modulo   inversion of the angles $\alpha$ and $\beta$.

Therefore, we can view (\ref{eq:spinrot}) directly (instead of its dual) as a non-commutative consensus iteration $\Phi$. We compute that $\Delta(\Phi^k) < \infty$ for $k \geq 2$, except if (a) $\alpha = 0 {\rm mod} \pi$ or (b) $\beta = 0 {\rm mod} \pi$, or (c) $(2 \alpha,\, 2\beta) = (0,0) {\rm mod} \pi$. Thus except for these special cases, we can apply Proposition 1 and Theorem 3 to conclude that $Z$ converges to $\tfrac{1}{2}\, I$. The physical interpretation is that, when stochastically alternating between two generic rotations, the probability distribution for the state becomes completely uniform after many iterations.

A Kraus map that is  both unital and trace-preserving can be considered as the analog of a doubly stochastic map in  classical consensus theory. The consensus value $1/2$ for a trace-perserving two-dimensional  quantum map is reminiscent of the   average-preserving property of doubly stochastic maps.

For the particular value $\beta = \frac{\pi}{2}$, the set of diagonal matrices $X$ is invariant under (\ref{eq:spinrot}). In fact,
$$\Phi(X) = \left(\begin{array}{cc} p x_{11} + (1-p) x_{22} & 0 \\ 0  & p x_{22} + (1-p) x_{11} \end{array} \right) $$
for a diagonal matrix $ X = \text{diag} \{ x_{11}, x_{22} \} $
which can be interpreted as a classical consensus algorithm  for the doubly-stochastic map
$$A = \left(\begin{array}{cc} p & (1-p) \\ (1-p)  & p \end{array} \right) \; .$$
$\hfill \square$\vspace{2mm}

\noindent \textbf{Example 3: spontaneous emission.} Consider a two-level quantum system as in Example 2. The expectation trajectory of the system undergoing the stochastic process of spontaneous emission is described by
\begin{eqnarray}
\label{eq:SpontEm}
\nonumber	
V_0 & := & V_{\text{no emission}} = \left( \begin{array}{cc} 1 & 0 \\ 0 & \sqrt{1-\gamma^2} \end{array} \right )
 , \quad V_1  :=  V_{\text{photon emission}} = \left( \begin{array}{cc} 0 & \gamma \\ 0 & 0 \end{array}\right)
\end{eqnarray}
with small $\gamma$ equal to the quotient of the chosen time discretization step by  the expected lifetime of the  unstable excited state. The dual of the Kraus map is
$$\Phi \; : \; X = \left(\begin{array}{cc} x_{11} & x_{12} \\ x_{12}^*  & x_{22} \end{array} \right) \; \longrightarrow
\Phi(X) = \left(\begin{array}{cc} x_{11} & x_{12} \sqrt{1-\gamma^2} \\ x_{12}^* \sqrt{1-\gamma^2} & x_{22}(1-\gamma^2) + x_{11} \, \gamma^2 \end{array} \right) \, .$$
From the simple form of the map, it is easy to guess asymptotic convergence to the consensus state $x_{11} I$. This property is reminiscent of the consensus property of a directed graph where the root agent eventually imposes its value to all others. Indeed, restricting $\Phi$ to diagonal $X$ results in the commutative consensus map of Example 1. (Following the construction proposed in Section VI.C, one has $V_0=S_0W_0$ and $V_1=S_1W_1$ with $W_0 = {\rm diag}\{0,\gamma\}$ and $W_1 = {\rm diag}\{1,\sqrt{1-\gamma^2}\}$.)

As in Example 1, this is a situation where the diameter of the quantum stochastic map is infinite for any finite iteration, yet the contraction of
spectrum range can be computed explicitly. The calculation gives
$$
[ \lambda_{\min}(\Phi(X)), \lambda_{\max}(\Phi(X))] =   [\lambda_{\min}(X) +   \rho_+,
 \lambda_{\max}(X) - \rho_-]
$$
for $ \rho_{\pm} = \gamma^2 \left( \sqrt{(\tfrac{x_{11}-x_{22}}{2})^2 + \vert x_{12} \vert^2} \pm (\tfrac{x_{11}-x_{22}}{2}) \right)$. This
shows the strict decay of the Lyapunov function $V$ at each iteration. $\hfill \square$\vspace{2mm}



\section{Conclusion}
We studied some implications of Birkhoff theorem for the convergence analysis of consensus algorithms in two different cones: the positive orthant of the euclidean space and the cone of hermitian positive definite matrices. In the positive orthant, the application of Birkhoff theorem to linear stochastic maps provides an interpretation of Tsitsiklis Lyapunov function as the Hilbert distance to consensus in log coordinates. This connection does not seem to have been made in the literature. In the cone of positive definite matrices, the application of Birkhoff theorem provides a framework to study the convergence of quantum stochastic maps. This framework suggests to study the contraction of the interval $[\lambda_{\min},\lambda_{\max} ]$ over finite horizons. To the authors' knowledge, this approach is novel. The analogy with the classical consensus problem suggests that the potential of this approach is in situations where the stochastic map is not self-adjoint (the analog of {\it directed} communication graph
in classical consensus theory) and/or in situations where the stochastic map is time-dependent (the analog of {\it time-varying} communication graphs in the classical consensus theory). These are the situations where common quadratic Lyapunov functions may fail to exist. The analog of Theorem 1 in the non-commutative setting remains to be characterized.


\section{Acknowledgments}

The first and second authors want to thank the Centre Automatique et Syst\`emes at Mines ParisTech for its hospitality and stimulating research environment when they were visiting in 2009-2010. This paper presents research results of the Belgian Network DYSCO (Dynamical Systems, Control, and Optimization), funded by the Interuniversity Attraction Poles Programme, initiated by the Belgian State, Science Policy Office. The scientific responsibility rests with its authors. A.Sarlette is supported as an FNRS postdoctoral researcher (Belgian Fund for Scientific Research). This work was also  supported in part by the "Agence Nationale de la Recherche" (ANR),
    Projet Blanc  CQUID number 06-3-13957.\\



\end{document}